\documentclass[10pt]{article}
\usepackage{amsfonts, amssymb, amsmath, amscd, longtable, tabularx, multicol}
\usepackage{graphicx}
\newtheorem{theorem-intro}{Theorem}[section]
\newtheorem{theorem}{Theorem}

\newtheorem{corollary-intro}[theorem-intro]{Corollary}
\newtheorem{remark}{Remark}

\newtheorem{example}{Example}

\def\C{\mathcal{C}^1(\Sigma)}
\def\G{\mathcal G}

\def\S{\Sigma}
\def\a{\alpha}
\def\b{\beta}
\def\g{\gamma}
\def\d{\delta}
\def\w{\omega}

\begin{document}
\title{The curve complex has dead ends}
\author{Joan S.Birman\footnote{supported in part by Simons Foundation Award \# 245711} \  and William W. Menasco}
\date{Geom. Dedicata DOI 10.1007/s10711-014-9978-y}
\maketitle

\begin{abstract} 
\noindent  It is proved that the curve graph $\C$ of a surface $\S_{g,n}$ has a local pathology that had not been identified as such: there are vertices $\a,\b\in\C$ such that $\b$ is a dead end of every geodesic joining $\a$ to $\b$.  There are also double dead-ends.   Every dead end has depth 1.  
\end{abstract}



\section{Introduction}\label{S:introduction}
The {\it curve graph} $\C$ of a surface $S$ is a graph whose vertices are homotopy classes of simple closed curves on $S$, with two vertices joined by an edge when the curves have disjoint representatives.  Assigning length 1 to each edge allows us to define the distance between two vertices as the length of a shortest path between two vertices.  The graph $\C$ is the 1-skeleton of the {\it curve complex} $\mathcal C(\Sigma)$,  a metric space and a simplicial complex that was introduced in the early 1970's by W. Harvey \cite{Ha}.  During the past 15 years a new and highly successful attack was made, in the groundbreaking papers of Masur and Minsky \cite{MM-I,MM-II}, in understanding the large-scale or `coarse'  geometry of $\mathcal C(\Sigma)$.   In that work it was necessary to overcome difficulties created by two types of  local pathology in $\C$.   The first is that $\C$ is locally infinite, that is there are infinitely many vertices that are distance 1 from any given vertex.  The second is that, typically, there are infinitely many distinct geodesics joining vertices $\a,\b$, where geodesics are regarded as being distinct when there is no homeomorphism of the surface taking the curves on one to the curves on the other.  
In this paper we will prove that there is yet another local pathology which, surprisingly, seems to have not been noticed as such.   

A vertex $\b\in\C$ is a {\it dead end} with respect to a vertex $\a$  if no geodesic $\G$ joining $\a$ to $\b$ can be extended past $\b$ to a longer geodesic.   
The {\it depth} of a dead-end in a geodesic of length $n$ is the number of edges one must follow backwards from the dead end to reach a vertex from which the shortened geodesic can be extended to a geodesic of length $\geq n+1.$   
These concepts are familiar ones in geometric group theory, where groups have been identified whose Cayley graphs have dead ends, and even dead ends of arbitrarily large depth, e.g. see \cite{Wa}.  
The main result in this paper is to prove that dead ends exist (and indeed are ubiquitous) in $\C$.  On the other hand, we will also prove that every dead end has depth 1.

To state our results, let $\Sigma =  \Sigma_{g,n}$ be an orientable  surface of genus $g$ with $n$ boundary components or punctures,  $3g-3+n >0$.   A curve on $\Sigma$ is {\it essential} if  does not bound either a disc or an annulus parallel to a component of $\partial\Sigma$.   In this note we will prove: 
\begin{theorem} \label{T:main theorem}   Let $\a,\b$ be essential curves on $\S$, and also (by an abuse of notation) vertices in $\C$.  
Then the following hold:
\begin{enumerate}
\item  [{\rm (A)}] If $\b$ is non-separating as a curve on $\S$, then $\b$ cannot be a dead-end with respect to any vertex $\a$ in the curve graph.  
\item [{\rm (B)}] Let $\a,\b\in\C$ with $n=d(\a,\b)\geq 3$.  Then necessary and sufficient conditions for $\b$ to be a dead end  with respect to $\a$ are {\rm (i)} the curve $\b$ separates $\Sigma$ in such a way that both components of $\S$ split along $\b$ support essential curves, and {\rm (ii)}
there exist distinct geodesics, $\G^{(1)}$ and $\G^{(2)}$ joining $\a$ to $\b$, such  that the vertices that immediately precede $\b$ on $\G^{(1)}$ and $\G^{(2)}$  are in distinct components of $\Sigma$ split along $\beta$.   
\item [{\rm (C)}] Let $\a$ be any non-separating curve on $\Sigma$.  Then  every geodesic of length $k\geq 0$ that ends at $\a$  can be extended to a geodesic of length $k+2$ that terminates in a dead end.  Indeed, there are infinitely many such extensions.   Moreover, if $\mu$ is also non-separating,  then any geodesic joining $\mu$ to $\a$  has infinitely many extensions  to a geodesic of length $d(\mu,\a)+4$ that has a double dead-end.   
\item [{\rm (D)}] Every dead end in $\C$ has depth 1.
 \end{enumerate}
\end{theorem}

\begin{remark} {\rm
We first encountered dead-ends when we were searching for examples of vertices $\a,\b\in\C$ with $d(\a,\b)=4$, and discovered that the only thing that was well-known was how to construct examples with $d(\a,\b)\geq 3$.   Attempting to extend a known distance 3 example to distance 4, we encountered a dead-end.   We asked several experts whether the existence of dead-ends was known, and to our surprise they all said no.  As we learned more, however, we realized that while the phenomenon had not been recognized as such, it had indeed affected many proofs.\footnote{One example is on page 910-912 of \cite{MM-II}, and in particular lines $8^--7^-, p. 910$, where we first  learned about the consequences of the fact that $\mathcal{C}(\S)$ is $\delta$-hyperbolic \cite{MM-I} for fellow-traveller geodesics.  We pondered over the statement `If $v$ is non-separating', wondering what would go wrong if $v$ was separating?  This case is handled separately on page 911.  Knowing about dead-ends could have made that discussion easier to understand intuitively.}
} \end{remark}  

\begin{example}
\label{E:basic example}  
{\rm We give a very simple family of examples that illustrate the existence of dead-ends and double dead-ends.   Let $\a,\b\subset\S_{g,n}$ be essential curves, with $\b$ separating and $\a\cap\b\not=\emptyset$, such that $\a\cup\b$ does not fill either component of $\S_{g,n}$ split along $\b$.   See the left sketch in Figure~\ref{F:basic example} for an example on $\S_{2,0}$.    We have chosen $\alpha=\mu$ for an example with $\alpha$ non-separating, and $\alpha=T_\mu(\beta)$, where $T_\mu$ means a Dehn twist along $\mu$,  for an example with $\alpha$ separating. 
Since $\a\cap\b\not=\emptyset$, but $\a\cup\b$ do not fill $\S_{g,n}$, the distance $d(\a,\b)\geq 2$.  Even more, since $\a\cup\b$ do not fill {\it either side} of $\S_{2,0}$ split along $\b$, we may find essential curves  $\d^{(1)}, \d^{(2)}$, one on each side of $\b$, giving paths $\G^{(1)} =(\a\to\d^{(1)}\to\b)$ and $\G^{(2)} = (\a\to\d^{(2)}\to\b)$.  Thus $d(\a,\b)=2$ and $\G^{(1)},\G^{(2)}$ are geodesics.   

Suppose that a step away from $\b$ to $\gamma$ could take us a step further from $\a$.  Then both $(\a\to\d^{(1)}\to\b\to\gamma)$ and $(\a\to\d^{(2)}\to\b\to\gamma)$ are length 3 geodesics.  Since $\g\cap\b=\emptyset$, $\g$ must be in $\S^{(1)}_{2,0}$ or $\S^{(2)}_{2,0}$, say $\S^{(1)}_{2,0}$.  But then $\d^{(2)}\cap\g=\emptyset$, which implies that $\a\to\d^{(2)}\to\g$ is a path, so that $d(\a,\g) =2$, not $3$.  So $\g$ must be in $\S^{(2)}$.  But then $\a\to\d^{(1)}\to\g$ would be a path of length 2 and again $d(a,\g)=2$.   Thus no such $\g$ exists, and $\b$ is a dead end with respect to $\a$.  

If we choose $\alpha=T_\mu(\beta)$, which is a  separating curve,  the identical argument, using the same two paths,  proves that both $\a$ and $\b$ are dead-ends. 

The right sketch in Figure~\ref{F:basic example} shows that if the same example is modified by adding handles and punctures, the proof will go through unaltered.}
\begin{figure}[htpb!]
\label{F:basic example}
\centerline{\includegraphics[scale=.58] {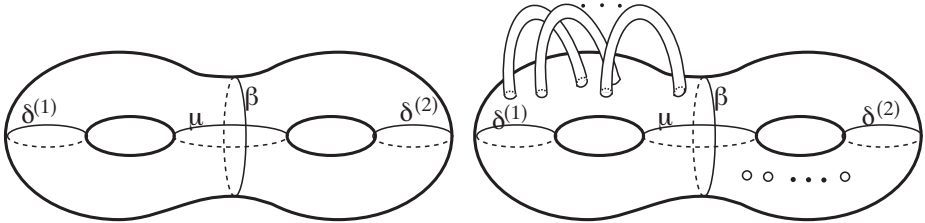}}
\caption{A family of length 2 geodesics with dead ends on $\S_{g,n}, g\geq 2, n\geq 0$.}
\end{figure} 
\end{example}

\section{Proof of Theorem~\ref{T:main theorem}}   We consider the proofs of (A),(B),(C) and (D) separately.

{\bf Proof of (A)}  The key observation that is needed for the proof is given in lines 8- and 7- on page 910 of \cite{MM-II}:  Assume that $\b$ is non-separating.  Choose a vertex $\g\in\C$ such that that $d_{\Sigma - \b}(\a,\g) >M$, where $d_{\Sigma - \b}$ means distance under subsurface projection and where
$M = M(\Sigma)$ is the constant  given in the Bounded Geodesic Image Theorem of \cite{MM-II}.  Then, as is proved on pages 910-911 of \cite{MM-II}, every geodesic  joining $\a$ to $\g$  must pass through $\b$, i.e.  there exists an extension of every geodesic joining $\a$ to $\b$.   By making different choices of $\g$ we obtain infinitely many such extensions for each fixed $\a$.   

{\bf Proof of (B):}    From the discussion regarding Example~\ref{E:basic example}  we know that (i) and (ii)  are sufficient for $\b$ to be a dead end with respect to $\a$. From the proof of part (A) of Theorem~\ref{T:main theorem} we also know that (i) is necessary.  It remains to prove  that (ii) is necessary, i.e. if $\b$ is separating and a dead end with respect to $\a$, then there always exist two geodesic paths $\G^{(1)}, \G^{(2)}$  with the given properties.    (This is the only place in the proof of (B) where we need $n=d(\a,\b)\geq 3$. )

Let $\Sigma^{(1)}, \Sigma^{(2)}$ be the 2 components of $\Sigma$ split along the dead end vertex $\b$. The fact that  $n\geq 3$  implies that $\a$ and $\b$ fill $\Sigma$, therefore $\a$ intersects both $\Sigma^{(1)}$ and $\Sigma^{(2)}$.    Assume that $\mathcal C(\Sigma^{(1)}), \mathcal C(\Sigma^{(2)})$  both have infinite diameter.  We may choose a non-separating curve $\w\subset\Sigma^{(2)}$ with $d_{\Sigma - \Sigma^{(1)}}(\a,\w)\geq 2n+1$.  Lemma 2.2 of \cite{Sch} asserts that, in this situation, every geodesic  joining $\alpha$ to $\omega$ has a vertex, say $\delta^{(1)}$, that (regarding that vertex as a curve on a surface) misses $\Sigma^{(2)}$, i.e. $\delta^{(1)}\subset\Sigma^{(1)}$.  Since $\w\subset \Sigma^{(2)}$, we have learned that $\d^{(1)}\cap \w= \d^{(1)}\cap \b = \emptyset$.  But then $d(\d^{(1)},\w) = d(\d^{(1)},\b) =1$.       

We wish to determine $d(\a,\d^{(1)})$.  Observe that since $\beta=\Sigma^{(1)}\cap\Sigma^{(2)}$, also $\w\subset \Sigma^{(2)}$ and $\d^{(1)}\subset\Sigma^{(1)}$, we know that any two of the 3 vertices $\b, \w, \d^{(1)}$ are distance 1 apart in $\C$.   By the triangle inequality,
$d(\a,\b) = n  \leq d(\a,\d^{(1)}) + 1.$  \  But then, $n-1  \leq d(\a,\d^{(1)})$, which implies that $n\leq d(\a,\w)$.  However, since $\b$ is a dead end, and $\b\cap \w = \emptyset$,  we also know that $d(\a,\w)\leq n$,  so  $d(\a,\w) = n$.  But then $d(\a, \d^{(1)})= n-1$.  
Choose any geodesic path joining $\a$ to  $\d^{(1)}$.  That path extends to the first sought-for  geodesic path $\G^{(1)} = (\a,\dots,\d^{(1)},\b)$.   
Now observe that if we had assumed that $\w\subset \Sigma^{(1)}$ instead of $\w\subset \Sigma^{(2)}$, we could have applied the same argument to obtain  a second geodesic path $\G^{(2)} = (\a,\dots,\d^{(2)},\b)$.  Since $\d^{(1)}$ and $\d^{(2)}$ are separated by $\b$, we have constructed the required two paths, and proved (B).

{\bf Proof of  (C):} We are given a geodesic $\G = (\g_0,\dots,\g_k)$, which ends in the non-separating vertex $\g_k$.   Let $\G^{(1)},\G^{(2)}$ be the length two geodesics that were constructed in Example~\ref{E:basic example}, where (since any two non-separating curves are equivalent under a homeomorphism of $\Sigma$) we may assume without loss of generality that the curve $\a$ in Example~\ref{E:basic example} is $\g_k$.   The concatenated  paths $\G\circ\G^{(1)}$ and $\G\circ \G^{(2)}$, both extend $\G$, but they might not be geodesics.  We now consider $(\Sigma\setminus \a)$, the surface $\Sigma$ split along the curve $\a=\g_k$ and the subsurface projection $\pi:\Sigma\to(\Sigma\setminus \a)$.  Choose a pseudo-Anosov map $f$ of $(\Sigma\setminus \a)$.  
For $q>M$ the concatenations of $\G$ with   $f^q( \G^{(1)})$ and $f^q(\G^{(2)})$  extend $\G$ to distinct geodesics of length $k+2$ that join $\g_0$ to $f^q(\b)$.  Since $\b$ is separating, we know that $f^q(\b)$ is too.  By (B), neither geodesic path has an extension beyond $f^q(\b)$. 

To build infinitely many such extensions of $\G$, choose distinct pseudo-Anosov maps  $f_1,f_2,,\dots$  and powers $q_1,q_2,\dots$ where each $q_i>M$. Then the geodesics $\G\circ f_i^{q_i}(\G^{(j)}), \ i=1,2,\dots, j=1,2$ give distinct extensions of $\G$ to geodesics of length $k+2$,  where the latter all have dead ends.   The construction can be done at each endpoint of $\G$ when both $\g_0$ and $\g_k$ are non-separating curves, to give double dead ends.      

{\bf Proof of (D):} Choose any geodesic $\G = (\g_0,\dots,\g_n)$ of length $n$ with a dead-end at $\g_n$.  We have just seen, in Example 1 and our proof of (C),  that the length $n-1$ subgeodesic $(\g_0\to\cdots\to\g_{n-1})\subset \G$  can be extended to a path of length $n+1$, with the middle curve, call it $\delta$, in the length 2 extension a non-separating curve that is distance $n$ from the initial curve $\g_0$.  By (A), the geodesic $(\g_0,\g_1,\dots,\g_{n-1},\d)$ of length $n$ can be extended arbitrarily far after that.  Since it includes the subpath $(\g_0,\dots,\g_{n-1})$, we have proved (D).   
{\hfill $\square$}\   


\noindent  {\bf Acknowledgements:}  We thank Dan Margalit, and  Jason Manning for many conversations about the Masur-Minsky machinery; also Yair Minsky for suggesting to us that part (A) of Theorem 1 ought to be true; also Saul Schleimer for suggesting how to simplify Example~\ref{E:basic example} from our original length 3 examples to ones of length 2.

\bigskip

\noindent
{\sc Joan Birman\\
Department of Mathematics, Barnard-Columbia\\
2990 Broadway\\
New York, NY 10027, USA\\}
\medskip
{\rm jb@math.columbia.edu}

\noindent
{\sc William W. Menasco\\
Department of Mathematics\\
University at Buffalo--SUNY\\
Buffalo, NY 14260-2900, USA\\}
\medskip
{\rm menasco@buffalo.edu}

\end{document}